\documentclass[11pt]{article}       
\usepackage{graphicx}
\usepackage{amsmath,amssymb,amsfonts,amsbsy}
\DeclareMathOperator\arcsinh{arcsinh}
\newtheorem{definition}{Definition}
\newtheorem{theorem}{Theorem}
\newtheorem{remark}{Remark}
\newenvironment{proof}{{\bf Proof.}}{\begin{flushright}$\blacksquare$\end{flushright}}
\newenvironment{example}[1]{\vspace{1ex}\noindent {\bf Example #1}}{}

\newcommand{\dword}[1]{\textsf{#1}}
\newcommand{\cs}{\mathcal{S}}
\newcommand{\TT}{\mbox{\boldmath $T$}}
\newcommand{\sig}{\mbox{\boldmath $\sigma$}}

\newcommand{\mm}{\boldsymbol{m}}
\newcommand{\es}{\bar{\sigma}}

\newcommand{\bL}{\bar{L}}
\newcommand{\hh}{\mbox{\boldmath $h$}}
\begin{document}

\title{Quasi-collocation based on CCC--Schoenberg operators and collocation methods
}


\author{Tina Bosner
}




\maketitle

\begin{abstract} 
We propose a collocation and quasi-collocation method for solving second order boundary value problems $L_2 y=f$, in which the differential operator $L_2$ can be represented in the product formulation, aiming mostly on singular and singularly perturbed boundary value problems. 
Seeking an approximating Canonical Complete Chebyshev spline $s$ by a collocation method leads to demand that $L_2s$ interpolates the function $f$.
On the other hand, in quasi-collocation method we require that $L_2 s$ is equal to an approximation of $f$ by the Schoenberg operator. 
We offer the calculation of both methods based on the Green's function, 
and give their error bounds.
\end{abstract}

\section{Introduction}
 
We are interested in solving boundary value problem
{\setlength\arraycolsep{2pt}
\begin{eqnarray}
L_2 y & = & f \ {\rm on}\ (a,b), \label{boundvalue} \\
y(a) & = & \alpha, \ y(b)=\beta, \nonumber
\end{eqnarray}}%
where $L_2$ is differential operator of order $2$, which can be represented in product formulation 
\begin{equation} \label{productdifop}
L_2 y:=p_1 D(p_2 D(p_3 y))
\end{equation}
for some $p_1,p_2,p_3$ of constant sign, 
with $D$ denoting ordinary derivative, like in \cite{bosnerroginanonuniformtension,bosnerroginasingular,bosnerad}. 

For example, each linear second order differential operator
$$
L_2y = y''+p_1 y'+p_2 y
$$
with non constant coefficients, can be represented in the form
$$
L_2 y=(py')'+q y,
$$
where $p(t):=e^{\int_a^tp_1(s)\,ds}$ and $q(t):=p_2(t)p(t)$. By \cite{coppeldisconjugacy}, there exists $\ell>0$ such that on each segment with length smaller than $\ell$, we can find an $u\in\ker{(L_2)}$ without zeros. In that case we have
\begin{equation} \label{productop}
L_2 y=u^{-1}D\bigg(p u^2 D\left(\frac{y}{u}\right)\bigg).
\end{equation}
If our segment $[a,b]$ is such that $b-a<\ell$, we can transform the linear equation in the product form.
Also, if the equation in (\ref{boundvalue}) is linear with constant coefficients, we can use an another approach to do the same thing, which says that there exist a \dword{critical length} $\ell$ (see \cite{carnicermainarpenacritical,carnicerpenacritlengthcycloid,beccaricasciolamazurecritical}) which determines that the kernel of the operator $L_2$ on some interval $[a,b]$ is an Extended Chebyshev space, if the length of $[a,b]$ is less than $\ell$. In that case, it is possible to transform the operator $L_2$ in the product representation like in (\ref{productdifop}) (see \cite{coppeldisconjugacy,mazurefinding}).
Specially, we are aiming here on singular and singularly perturbed boundary value problems.

One way of solving such a problem are the collocation methods \cite{russell1,deBoor:1,ahlbergito,kadalbajofitted,marusicrogina,marusiccollocation,bosnerroginasingular,bosnerad}, in which we are searching for $s\in\cs$, where $\cs$ is, usually, some spline space, such that
$$
L_2s(\tau_i)=f(\tau_i), \quad i=1,\dots,n,
$$
for some choice of collocation points $\tau_i$, $i=1,\dots,n$, $n:=\dim{(\cs)}$. In the collocation method, $L_2 s$ is actually interpolated, so next to the interpolation, the idea is to replace it, with some other kind of approximation, to be precise, with the approximation by CCC--Schoenberg operators \cite{bosnerroginaschoenberg} (see also \cite{mazurechebschoenberg}). Something similar for the polynomial splines was done in \cite{FOUCHER20093455}. 

The paper is divided in five sections. In Section \ref{sec1} we gathered all facts and theorems we need to define our methods: starting from defining CCC--space and CCC--splines, followed by CCC--Schoenberg operator, and finally interpolation. In the next section we introduce the collocation and quasi-collocation method and give their error bounds. We also suggest two ways of calculating these approximations, and then test  them in Section \ref{secexample}. 
We end with the conclusion.
 
\section{Preliminaries} \label{sec1}
The main tool for our method are splines. To define a spline space, at the beginning we have to determine where the pieces of a spline are coming from. 

\subsection{CCC--space}
First, we specify the space that we consider. For more details and proofs check \cite{bosnerroginaschoenberg,bosnerphd}.
\begin{definition} \label{defCCC}
	Let $u_1$ be a bounded positive function on the interval $[a,b]$, and let there be
	given a \dword{measure vector} $d\sig:=(d\sigma_2,\ldots
	d\sigma_{k})^{\rm T}$, where $\sigma_2,\dots,\sigma_k$ are bounded, continuous and strictly increasing functions on $[a,b]$.  Let
	{\setlength\arraycolsep{2pt}
		\begin{eqnarray}
			u_2(x) & = & u_1(x)\int_{a}^x d\sigma_{2}(\tau_{2}), \nonumber \\
			& \vdots & \label{CCCs} \\
			u_{k}(x) & = & u_1(x)\int_{a}^x d\sigma_{2}(\tau_{2}) \int_{a}^{\tau_2} \ldots \int_{a}^{\tau_{k-1}}
			d\sigma_{k}(\tau_{k}), \nonumber
		\end{eqnarray}}%
		for $x\in[a,b]$. We call $U_k=\{u_i\}_1^k$ a \dword{Canonical Complete Chebyshev (CCC)--system}.
\end{definition}
In the most general case, it is enough to assume that  $\sigma_i$ are only right continuous~\cite{schumaker1}. However, in the
sequel we shall use measures that posses densities with respect to the Lebesgue measure (possibly not smooth, or not even continuous), and henceforth assume that
$\sigma_i$ are continuous. As a consequence, every $u_i$, $i=1,\dots,k$ from Definition \ref{defCCC} is continuous. These densities do not have to be strictly positive: they can posses some final number of zeros, and they can even have some poles (see for examples \cite{bosnerroginavariable,bosnerroginasingular}), as long as their distribution function are bounded and strictly increasing.
If for the density $w_i$ of the measure $d\sigma_i$, we have $w_i>0$ and $w_i\in C^{k-i}([a,b])$, and this is true for each $i=2,\dots,k$, then this CCC--system is also an \dword{Extended Complete Chebyshev (ECC)--system} \cite{schumaker1}.

%
%

The $i^{\rm th}$ \dword{reduced system} is defined to be a Chebyshev system
corresponding to the \dword{reduced measure vector}, that is
$$
d\sig^{(i)} := (d\sigma_{i+2}, \ldots, d\sigma_k)^
{\rm T} 
,\quad i=1,\ldots k-2.
$$
We write $U_k^{(i)}:=\{u_{i,j}\}_{j=1}^{k-i}$ for $ i=1,2, \dots ,k-1$, where
{\setlength\arraycolsep{2pt}
\begin{eqnarray*}
	u_{i,1}(x) & = & 1, \\
	u_{i,2}(x) & = & \int_{a}^x d\sigma_{i+2}(\tau_{i+2}), \\
	& \vdots & \\
	u_{i,k-i}(x) & = & \int_a^x d\sigma_{i+2}(\tau_{i+2})\ldots \int_a^{\tau_{k-1}}
	d\sigma_{k}(\tau_{k}).
\end{eqnarray*}}%
%
\begin{definition} \label{CCCspace}
	The linear space
	$\cs(k,d\sig,u_1):={\rm span}\{U_k\}$ is called a \dword{CCC--space} if $U_k$ is a CCC--system.
	If $u_1\equiv 1$ then $\cs(k,d\sig):=\cs(k,d\sig,u_1)$.
\end{definition}
Also, if $U_k$ is an ECC--system, then the space form Definition \ref{CCCspace} is called an \dword{ECC--space}.
\begin{definition} \label{genderiv}
	Assuming that each of the functions $\sigma_2$, \dots ,$\sigma_k$ have been extended
	to an interval to the left of $a$ and to the right of $b$, \textit{i.e}. to an interval $[c,d]\supset
	[a,b]$ with $c<a$, $d>b$, 
	then for
	$f\in\cs(k,d\sig,u_1)$ and $x\in[a,b]$ we define \dword{generalized derivatives} as
	$$
	L_j:=D_{j}\cdots D_{1}D_0, \quad j=0,\ldots ,k-1,
	$$
	where
	{\setlength\arraycolsep{2pt}
	\begin{eqnarray}
	D_{0}f(x) & := & \frac{f(x)}{u_1(x)}, \nonumber \\
	D_{j}f(x) & := & \lim_{\delta \rightarrow 0^+}{\frac{f(x+\delta)-f(x)} 
		{\sigma_{j+1}(x+\delta)-\sigma_{j+1}(x)}},\ j=1,\ldots ,k-1. \label{defDj}
	\end{eqnarray}}%
\end{definition}

If each $\sigma_i$ posses a density, let us say, $w_i$, {\it i.e.} $d\sigma_i=w_i\,d\lambda$, then 
\begin{equation} \label{ECCder}
D_i=\frac{D}{w_{i+1}}, \quad i=2,\dots,n,
\end{equation} 
where $D$ stands for the ordinary derivative.
Also, as we mentioned before, $\frac{1}{w_{i+1}}$ does not have to be strictly positive or continuous, and as a consequence we allow functions $p_1$ and $p_2$ in (\ref{productdifop}) to have some zeros or poles.

For an arbitrary measure $d\es_{k+1}$, we can also define
$$
D_kf(x) := \lim_{\delta \rightarrow 0^+}{\frac{f(x+\delta)-f(x)}
	{\es_{k+1}(x+\delta)-\es_{k+1}(x)}},
$$
for $x\in[a,b]$, and $L_{k}:=D_{k}\cdots D_{1}D_{0}$. Now it is obvious that $\ker{(L_k)}=\cs(k,d\sig,u_1)$.
Then we have
$$ 
L_ju_i=\left \{ \begin{array}{ll}
0, & i=1,2, \dots ,j \\
u_{j,i-j}, & i=j+1, \dots ,k,
\end{array} \right. \qquad j=0,1, \dots ,k.
$$ 

The generalized derivatives can also be calculated for a function from a more general space of functions, like $C[a,b]$,
for which the required limits of the form (\ref{defDj}) exist.

%
%
%
\begin{theorem}
	CCC--space $\cs(k,d\sig,u_1)$ is $k$--dimensional linear space, with basis $\{u_1, u_2,$ $ \dots ,u_k\}$, and $L_{j}$
	is linear operator
	$\cs(k,d\sig,u_1) \rightarrow  \cs(k-j,d\sig^{(j)}):={\rm span}\{U_k^{(j)}\}$ for $j=1, \dots
	,k-1$.
\end{theorem}

\subsection{CCC--splines}
For a partition
$\Delta= \{x_{i}\}_{i=0}^{d+1}$ of an interval $[a,b]$,
given  \dword{multiplicity vector} $\mm=(m_{1},\ldots,$ $m_{d})$,
($0 < m_i \leqslant k$), and $M:=\sum_{i=1}^{d} m_i$, we shall denote by
$\TT:=\{t_{1} \ldots t_{2k+M}\}$ an \dword{extended partition}
in the usual way:
$$
t_1 = \cdots = t_k =  x_0=a,
$$
\begin{equation} \label{extendpart}
(t_{k+1},\dots,t_{k+M}) = (x_1^{[m_1]},\dots,x_d^{[m_d]}),
\end{equation}
$$
b=x_{d+1} = t_{k+M+1} = \cdots = t_{2k+M},
$$
where $x^{[m]}$ stands for $x$ repeated $m$ times.
\begin{definition}
	Let $\cs(k,d\sig,u_1)$ be a CCC--space. A set $\cs(k, d\sig,u_1, \TT)$ of bound\-ed
	functions on $[a,b]$ such that:
	\begin{itemize}
		\item[(i)] for each $s\in\cs(k, d\sig,u_1, \TT)$ and each $i=0,
		\dots ,d$, there exists $s_i\in\cs(k,d\sig,u_1)$ such that
		$s|_{[x_i,x_{i+1}]}=s_i|_{[x_i,x_{i+1}]}$,
		\item[(ii)]
		$L_js_{i-1}(x_i)=L_js_i(x_i)$ for $j=0, \dots ,k-m_i-1$, $i=1, \dots ,d$,
	\end{itemize}
	is called the \dword{space of Canonical Complete Chebyshev (CCC)--splines on extended partition $\TT$}. If $u_1\equiv 1$, then $\cs(k, d\sig, \TT):=\cs(k, d\sig,u_1, \TT)$.
\end{definition}
It is clear that $\cs(k, d\sig,u_1, \TT)$ is a linear space, and its dimension is equal to
$$
n:=\dim{\cs(k, d\sig,u_1, \TT)}=k+M.
$$
For the most convenient basis of the CCC--spline space $\cs(k, d\sig,u_1, \TT)$ we use the
\dword{B-splines} $T_j^k$, $j=1,\dots,n$ 
(see~\cite{bosnerroginaschoenberg,schumaker2,schumaker1}).

The CCC--B-splines posses the usual properties, listed in the following theorem.
\begin{theorem} \label{bsplprop}
	Let $\TT=\{t_i\}_{i=1}^{n+k}$ be an extended partition of $[a,b]$, and let $t_i<t_{i+k}$, $i=1,2,\dots,n$, then
	{\setlength\arraycolsep{2pt}
		\begin{eqnarray*}
		T_i^k(x) & = & 0, \quad x<t_i\ {\rm and}\ x>t_{i+k},\ i=1,2,\dots,n, \nonumber \\
		T_i^k(x) & > & 0, \quad t_i<x<t_{i+k},\ i=1,2,\dots,n, \nonumber \\
		\sum_{i=1}^n T_i^k(x) & = & u_1, \quad x\in[a,b]. \label{partitonunity}
		\end{eqnarray*}}%
\end{theorem}
%
%
\begin{remark} \label{partofunity}
Direct consequence of Theorem \ref{bsplprop} is that for $u_1\equiv 1$ the B-splines make partition of unity:
$$
\sum_{i=1}^n T_i^k(x) = 1
$$
for all $x\in[a,b]$.
\end{remark}

Let the B-splines associated with the $i^{\rm th}$
reduced space be denoted as
$T_j^{k-i}$, $i=0,\cdots,k-1$, and let
\begin{equation} \label{integralbspl}
C_j^{k-i}:=\int_{t_j}^{t_{j+k-i}} T_j^{k-i}(t)d\sigma_{i+1}(t).
\end{equation}
It is known that the de\thinspace{}Boor--Cox type recurrence does not exist for the general
Chebyshev splines~\cite{schumaker4}, but in many cases we can use the \dword{derivative formula} (\ref{ederivative})
instead~\cite{roginasingular,roginaderivative,bosnerphd}.
\begin{theorem}[Derivative formula] \label{derivativethm}
	Let $L_{1}$ be the first generalized derivative with respect to CCC--space
	$\cs(k,d\sig,u_1)$, and let the multiplicity vector $\mm=(m_1,\dots ,m_l)$ satisfy
	$m_i\leqslant k$ for $i=1, \dots, l$. Then for  $x\in [a,b]$ and $i=1, \dots, n$,
	the following derivative formula holds:
	\begin{equation}
	L_1T_i^{k}(x) =  \frac{T_i^{k-1}(x)}
	{C_i^{k-1}}-\frac{T_{i+1}^{k-1}(x)}{C_{i+1}^{k-1}}.
	\label{ederivative}
	\end{equation}
\end{theorem}
The formula can also be used to inductively define B-splines~\cite{bister}, although for stable calculation of CCC--splines,
the algorithms based on knot insertion are recommended
\cite{bosnerphd,roginafourth,roginabosnerlower,bosnerweighted,bosnerroginanonuniformtension,bosnerad,bosnerroginasingular,bosnerroginavariable}.


The difference between  the spaces $\cs(k, d\sig,u_1, \TT)$ and $\cs(k, d\sig, \TT)$ is that the functions from the first one are equal to the functions from the latter multiplied by $u_1$. The same is true for their B-splines.

\begin{theorem} \label{u1jed1}
The B-splines $T_i^k\in\cs(k, d\sig,u_1, \TT)$ , $\bar{T}_i^k\in\cs(k, d\sig, \TT)$ satisfy
$$
T_i^k(x) = u_1(x)\cdot\bar{T}_i^k(x), \qquad i=1,\dots,n,
$$
for all $x\in[a,b]$.
\end{theorem}
%
%
\begin{remark} \label{bsompu1}
	From the perspective of Theorem \ref{u1jed1}, from now on we can assume, without loss of generality, that $u_1\equiv 1$.
\end{remark}
\subsection{Marsden's identity} \label{sec2}
In the polynomial spline theory, Greville abscis\ae, or knot averages, $\{t_i^*\}_{i=1}^n$
are known as de\thinspace{}Boor points of the identity function, since
$$ 
x=\sum_{i=1}^n t_i^*B_i^k(x), \qquad t_i^*:=\frac{t_{i+1}+\cdots+t_{i+k-1}}{k-1},
$$ 
where $\{B_i^k\}_{i=1}^n$ are polynomial B-splines of order $k$ associated with the extended partition
$\{t_i\}_{i=1}^{k+n}$. 
In the Chebyshev setting, the function $u_2$ takes over the role of the identity function, and it can be shown that 
de\thinspace{}Boor points of $u_2$ can be obtained as the convex combination of $u_2(t_{i+1}),\dots,u_2(t_{i+k-1})$. By Theorem 18 from \cite{bosnerroginaschoenberg} we know how to represent each of the functions $u_j$ as a linear combination of  the B-splines
$$ 
	u_j(x)=\sum_{i=1}^{n} \eta_{i,j}^k T_i^k(x), \qquad j=1,2,\dots,k,
$$ 
and in particular, we know how to calculate $\eta_{i,2}^k$. Since the function $u_2$ is strictly increasing, $(\eta_{i,2}^k)_{i=1}^n$ is strictly increasing, as well as their originals, {\it i.e} the \dword{CCC--Greville points}
$(\xi_{i,2}^k)_{i=1}^n$ ($\eta_{i,2}^k=u_j(\xi_{i,2}^k)$). Even more, we have that $\xi_{i,j}^k\in[t_{i+1},t_{i+k-1}]$.

\subsection{CCC--Schoenberg operators}

In the polynomial case, Greville abscis\ae{}  
define classical Schoenberg operator applied to a function $f\in C[a,b]$
$$
S[f]:=\sum_{i=1}^n f(t_i^*)B_i^k,
$$
where $\{B_i^k\}_{i=1}^k$ are as before.
When working with CCC--splines, as we mentioned, the function $u_2$ takes over the role of the identity function, and $(\xi_{i,2}^k)_{i=1}^n$ of the Greville abscis\ae{}. But, the generalization of the Schoenberg operator goes in two directions: first, by replacing the polynomial setting by CCC one, and second, the Greville abscis\ae{} can be replaced by a more general set of points than just $(\xi_{i,2}^k)_{i=1}^n$. For proofs and details see \cite{bosnerroginaschoenberg}.


\begin{definition} \label{Schoenbergdef}
	A linear operator $S$ on $C[a,b]$ is called \dword{a CCC--Schoenberg operator based on the CCC--spline space
		$\cs(k,d\sig,\TT)$} if
	\begin{itemize}
		\item[(i)] it is of the form
		\begin{equation} \label{defeqshoenberg}
		S[f]:=\sum_{i=1}^n f(\zeta_i)T_i^k
		\end{equation}
		for all $f\in C[a,b]$, where the nodes $\zeta_1,\dots,\zeta_n$ satisfy $a=\zeta_1<\zeta_2<\cdots<\zeta_{n-1}<\zeta_n=b$,
		\item[(ii)] it reproduces a two-dimensional CCC--space $\cs_2\subset\cs(k,d\sig,\TT)$ in the sense
		$$
		S[f]=f,\quad {\rm for}\ \forall f\in\cs_2.
		$$
	\end{itemize}
\end{definition}

		Because the B-splines satisfy the partition of unity, the operator $S$ defined by (\ref{defeqshoenberg}) must reproduce
		the constants. Also, since $\cs_2$ is a CCC--space, there exists a continuous, strictly increasing distribution function $\es_2$ such that
		$\cs_2=\cs(2,d\bar{\sig}_2)$, with $d\bar{\sig}_2:=(d\es_2)^{\rm T}$. That means that actually $\cs_2={\rm span}\{1,\es_2\}\subset
		\cs(k,d\sig,\TT)$, and $\es_2=\sum_{i=1}^n \nu_iT_i^k$ for some $\nu_1,\dots,\nu_n$.

The following theorem determines which strictly increasing splines can be reproduced by a CCC--Schoenberg operator.
\begin{theorem} \label{strictdeboor}
	For a given continuous spline $s:=\sum_{i=1}^n s_iT_i^k\in\cs(k,d\sig,\TT)$, the following two properties are equivalent:
	\begin{itemize}
		\item[(i)] the spline $s$ is strictly increasing on $[a,b]$ and it is reproduced by a unique CCC--Schoenberg
		operator based on $\cs(k,d\sig,\TT)$;
		\item[(ii)] the de\thinspace{}Boor points $s_i$, $i=1,\dots,n$, of $s$ form a strictly increasing sequence.
	\end{itemize}
\end{theorem}

\begin{remark}
\begin{itemize}
\item[(i)] The connection between $(s_i)_{i=1}^n$ from Theorem \ref{strictdeboor} and $(\zeta_i)_{i=1}^n$ from (\ref{defeqshoenberg}) is simple:
$$
s(\zeta_i)=s_i, \qquad i=1,\dots,n.
$$
\item[(ii)] It can be shown that $\zeta_i\in[t_{i+1},t_{i+k-1}]$.
\end{itemize}
\end{remark}

Let $\cs_2={\rm span}\{1,\es_2\}$, where $\es_2\in\cs(k,d\sig,\TT)$ satisfies the property (ii) of Theorem
\ref{strictdeboor}, $\es_2\in C[a,b]$, and let the CCC--Schoenberg operator
$$
S[f]:=\sum_{i=1}^n f(\zeta_i)T_i^k
$$
reproduces the CCC--space $\cs_2$. The question that arises is how well $S[f]$ approximates a function $f$
satisfying some mild conditions.

Let us for $f\in C[a,b]$ define
$$
\bL_1f(x):=\lim_{\delta\to 0^{+}}\frac{f(x+\delta)-f(x)}{\es_2(x+\delta)-\es_2(x)},
$$
and similarly as for dual derivatives 
(see the (3) from \cite{bosnerroginaschoenberg}), derivative
$$
\bL_1^- f(x):=\lim_{\delta\to 0^{+}}\frac{f(x)-f(x-\delta)}{\es_2(x)-\es_2(x-\delta)}
$$
as well as
$$
\bar{L}_2:=\bar{D}_2\bL_1, \qquad \bar{L}_2^-:=\bar{D}_2^{-}\bL_1^-,
$$
where
$$
\bar{D}_2 f(x):=\lim_{\delta\to 0^{+}}\frac{f(x+\delta)-f(x)}{\bar{\sigma}_3(x+\delta)-\bar{\sigma}_3(x)}, \qquad
\bar{D}_2^- f(x):=\lim_{\delta\to 0^{+}}\frac{f(x)-f(x-\delta)}{\bar{\sigma}_3(x)-\bar{\sigma}_3(x-\delta)},
$$
for an arbitrary $d\bar{\sigma}_3$. 
Finally, let $\Lambda_{(d\es_2,d\es_3)}([a,b])$ be the space of continuous
function on $[a,b]$ such that
\begin{itemize}
	\item[(i)] $\bL_1f(x)=\bL_1^-f(x)$, for $x\in[a,b]$,
	\item[(ii)] $\bar{L}_2f$ and $\bar{L}_2^-f$ are bounded and integrable on $[a,b]$ with respect to the $d\bar{\sigma}_3$.
\end{itemize}
If we have (\ref{ECCder}), and if the distribution function $\bar{\sigma}_2$ is smooth, as we mentioned in \cite{bosnerroginaschoenberg}, it should not be a big problem to check whether $f\in\Lambda_{(d\es_2,d\es_3)}([a,b])$.

Let
$$
\|f\|_\infty:=\max_{a\leqslant x\leqslant b}{|f(x)|},
$$
then for the approximation by CCC--Schoenberg operator to a function we can state the
following theorem.
\begin{theorem} \label{thmvdquadr}
	Let $f\in\Lambda_{(d\es_2,d\es_3)}([a,b])$, where $\es_2\in\cs(k,d\sig,\TT)$ is continuous and satisfies (ii) of Theorem
	\ref{strictdeboor}, while
	$d\bar{\sigma}_3$ is arbitrary. Let $S$ be the CCC--Schoenberg operator which reproduces $\cs_2={\rm span}\{1,\es_2\}$.
	Then
	$$ 
	\|f-S[f]\|_\infty\leqslant C \bar{h}_{(d\es_2,d\bar{\sigma}_3,\TT)}^2
	$$ 
	where
	$$ 
	\bar{h}_{(d\es_2,d\bar{\sigma}_3,\TT)}:=\max_{0\leqslant i\leqslant d}\max{\{\es_2(x_{i+1})-\es_2(x_i),\bar{\sigma}_3
		(x_{i+1})-\bar{\sigma}_3(x_i)\}},
	$$ 
	and $C$ does not depend on $\TT$.
\end{theorem}

\subsection{Interpolation}

Although, the interpolation problem in Chebyshev setting is well known, like for an example in \cite{mulbachinterpol,schumaker1,kochlycheinterpol,kayumovmazure,kvasovkimweight,CARNICER2014interpol,mazureinterpvsdesign}, we will pay some more attention to it, especially to the error bounds.
We follow the ideas from \cite{deboor1}. At the beginning, we will assume that all multiplicities of the interior knots are less then $k$ ($m_i<k$, $i=1,\dots,d$), so that $\cs(k,d\sig,\TT)\subset C([a,b])$.

We observe the following interpolation problem. For the function $f\in C([a,b])$, and the given interpolating points $(\tau_i)_{i=1}^n$, we search for a spline
$$
s(x)=\sum_{i=1}^n c_iT_i^k(x),
$$
$s\in\cs(k,d\sig,\TT)$, such that $s(\tau_i)=f(\tau_i)$. The \dword{interpolation operator} $I$ is such that $I[f]:=s$.

Let us denote with $A$ the interpolation matrix $A:=[T_j^k(\tau_i)]_{i,j}$. The same proof of the Schoenberg--Whitney Theorem for the polynomial case, works also for any CCC--spline space, so the matrix $A$ is regular if $T_i^k(\tau_i)\ne0$, $i=1,\dots,n$. Also, because of the support property of the CCC--B-splines, the interpolation matrix is banded, and because of Theorem 16 and 17 from \cite{bosnerroginaschoenberg}, $A$ is totally positive.

Now, we would like to give some error bounds of the interpolating spline. Let $f\in C([a,b])$, and we start with the fact that $I[g]=g$ for any $g\in\cs(k,d\sig,\TT)$, and therefore
{\setlength\arraycolsep{2pt}
\begin{eqnarray*}
\|f-I[f]\|_\infty & \leqslant & \|f-g\|_\infty+\|I[f-g]\|_\infty\leqslant\|f-g\|_\infty+\|I\|_\infty\|f-g\|_\infty \nonumber \\
 & = & (1+\|I\|_\infty)\|f-g\|_\infty. 
\end{eqnarray*}}%
By choosing the closest spline from $\cs(k,d\sig,\TT)$ to function $f$, we get
\begin{equation} \label{interpolgrocj1}
\|f-I[f]\|_\infty \leqslant (1+\|I\|_\infty)\,{\rm dist}\,(f,\cs(k,d\sig,\TT)).
\end{equation}

The calculation of bounds for ${\rm dist}\,(f,\cs(k,d\sig,\TT))$ follows mostly the derivation from \cite{deboor1}, with one exception, so we will concentrate more on this part. 
Let the operator $B$ be defined with
$$
B[g]:=\sum_{i=1}^n g(\tau_i)T_i^k,
$$
for some $\tau_1\leqslant\cdots\leqslant\tau_n$ and any continuous function $g$. This operator reproduces constants. Then for $x\in[t_j,t_{j+1}\rangle$ 
we have
{\setlength\arraycolsep{2pt}
\begin{eqnarray*}
|f(x)-B[f](x)| & \leqslant & \sum_{i=j+1-k}^j |f(x)-f(\tau_i)|T_i^k(x) \\
 & \leqslant & \max_{j-k+1\leqslant i\leqslant j}{\{|f(x)-f(\tau_i)|\}},
\end{eqnarray*}}%
then, by choosing the same $(\tau_i)_{i=1}^n$ as in \cite{deboor1}, 
we get, as before,
{\setlength\arraycolsep{2pt}
\begin{eqnarray}
\max_{j-k+1\leqslant i\leqslant j}\!\! & \{ & \!\!\!|f(x) - f(\tau_i)|\} \nonumber \\
 & \leqslant & \max\{|f(y)-f(z)|: y,z\in[t_{j+1-\frac{k}{2}},t_{j+1}]\ {\rm or} \nonumber \\
 & & y,z\in[t_j,t_{j+\frac{k}{2}}]\} \label{disttoder} \\
 & \leqslant & \left\lfloor \frac{k+1}{2} \right\rfloor\omega(f,\hh_1), \nonumber
\end{eqnarray}}%
where 
$$
\hh_1:=\max_{0\leqslant i\leqslant d}{\{x_{i+1}-x_i\}},
$$
and $\omega$ is modulus of continuity. Therefore it is again
\begin{equation} \label{distmodul1}
{\rm dist}\,(f,\cs(k,d\sig,\TT))\leqslant{\rm const}\, \omega(f,\hh_1).
\end{equation}

But, if the function $f$ has piecewise continuous $L_1f$, then we proceed differently. By the generalized Taylor expansion from Theorem 8 \cite{bosnerroginaschoenberg}, and the Mean value theorem for Lebesgue--Stieljes integrals (see \cite{burrill}) which says that there exists $d_{y,z}$, with $\min_{[y,z]}{|L_1f|}\leqslant d_{y,z}\leqslant\max_{[y,z]}{|L_1f|}$, 
such that from (\ref{disttoder}) it follows
{\setlength\arraycolsep{2pt}
\begin{eqnarray*}
\max_{j-k+1\leqslant i\leqslant j}\!\! & \{ & \!\!\!|f(x) - f(\tau_i)|\} \nonumber \\
 & = & \max{\left\{\left|\int_y^z L_1f(\tau_2)d\sigma_2(\tau_2)\right| :y,z\in[t_{j+1-\frac{k}{2}},t_{j+1}]\ {\rm or}\ y,z\in[t_j,t_{j+\frac{k}{2}}]\right\}} \nonumber \\
 & \leqslant & \max{\left\{\left|\int_y^z |L_1f(\tau_2)|d\sigma_2(\tau_2)\right| :y,z\in[t_{j+1-\frac{k}{2}},t_{j+1}]\ {\rm or}\ y,z\in[t_j,t_{j+\frac{k}{2}}]\right\}} \nonumber \\
 & = & \max{\left\{d_{y,z}\left|\int_y^z d\sigma_2(\tau_2)\right| :y,z\in[t_{j+1-\frac{k}{2}},t_{j+1}]\ {\rm or}\ y,z\in[t_j,t_{j+\frac{k}{2}}]\right\}} \nonumber \\
 & \leqslant & \|L_1f\|_\infty 
 \max{\left\{\left|\int_y^z d\sigma_2(\tau_2)\right| :y,z\in[t_{j+1-\frac{k}{2}},t_{j+1}]\ {\rm or}\ y,z\in[t_j,t_{j+\frac{k}{2}}]\right\}} \nonumber \\
 & = & \|L_1f\|_\infty \max{\{\sigma_2(t_{j+1})-\sigma_2(t_{j+1-\frac{k}{2}}),\sigma_2(t_{j+\frac{k}{2}})-\sigma_2(t_j)\}} \nonumber \\
 & \leqslant & \|L_1f\|_\infty \left\lfloor \frac{k+1}{2} \right\rfloor\hh_2, 
\end{eqnarray*}}%
since $\sigma_2$ is strictly increasing, and with
$$
\hh_2:=\max_{0\leqslant i\leqslant d}{\{\sigma_2(x_{i+1})-\sigma_2(x_i)\}}.
$$
All together, we get
\begin{equation} \label{nejednl1g}
{\rm dist}\,(f,\cs(k,d\sig,\TT))\leqslant{\rm const}\, \hh_2 \|L_1f\|_\infty.
\end{equation}

Now,
$$
{\rm dist}\,(f,\cs(k,d\sig,\TT))={\rm dist}\,(f-s,\cs(k,d\sig,\TT)), \quad \forall s\in\cs(k,d\sig,\TT),
$$
and from (\ref{distmodul1}) we get
$$
{\rm dist}\,(f,\cs(k,d\sig,\TT))\leqslant{\rm const}\, \omega(f-s,\hh_1), \quad \forall s\in\cs(k,d\sig,\TT).
$$
If $f$ and $s$ have piecewise continuous first generalized derivatives, then from (\ref{nejednl1g}) follows
$$
{\rm dist}\,(f,\cs(k,d\sig,\TT))\leqslant{\rm const}\,\hh_2\|L_1f-L_1s\|,
$$ 
for all $s\in\cs(k,d\sig,\TT)$. Again, we can choose $L_1s$ to be closest to the $L_1f$ in $\cs(k-1,d\sig^{(1)},\TT)$, to obtain the estimate
$$
{\rm dist}\,(f,\cs(k,d\sig,\TT))\leqslant{\rm const}\,\hh_2\,{\rm dist}\,(L_1f,\cs(k-1,d\sig^{(1)},\TT)).
$$
In case $L_1f$ is continuous, we could use (\ref{distmodul1})  again to get
$$
{\rm dist}\,(f,\cs(k,d\sig,\TT))\leqslant{\rm const}\,\hh_2\omega(L_1f,\hh_1).
$$
If $f$ has even continuous second generalized derivative, and if $k>2$, we can proceed in the same way to achieve
$$
{\rm dist}\,(f,\cs(k,d\sig,\TT))\leqslant{\rm const}\,\hh_2\hh_3\omega(L_2f,\hh_1),
$$
with generally
$$
\hh_j:=\max_{0\leqslant i\leqslant d}{\{\sigma_j(x_{i+1})-\sigma_j(x_i)\}}, \quad j=2,\dots,k.
$$
Let us just remark that all of the ``const'' in previous inequalities are not necessarily the same. 
That way we have proven the Jackson type theorem for CCC--splines.
\begin{theorem} \label{tmjackson}
For $j=0,\dots,k-1$, there exists constant {\rm const}, which depends on $k$ and $j$, such that for all $f$ with continuous $L_if$, $i=0,\dots,j$ on $[a,b]$
$$
{\rm dist}\,(f,\cs(k,d\sig,\TT))\leqslant{\rm const}\,\hh_2\cdots\hh_{j+1}\omega(L_jf,\hh_1),
$$
and for $j=k-1$ if $f$ has also continuous $L_kf$ for some arbitrary measure $d\es_{k+1}$, then 
$$
{\rm dist}\,(f,\cs(k,d\sig,\TT))\leqslant{\rm const}\hh_2\cdots\hh_k\bar{\hh}_{k+1}\|L_kf\|_\infty,
$$
with
$$
\bar{\hh}_{k+1}:=\max_{0\leqslant i\leqslant d}{\{\es_{k+1}(x_{i+1})-\es_{k+1}(x_i)\}}.
$$
\end{theorem}


Now we go back to (\ref{interpolgrocj1}), where
we can also say a little bit more about 
$$
\|I\|_\infty:=\max{\left\{ \frac{\|I[g]\|_\infty}{\|g\|_\infty}:g\in C([a,b])\backslash \{0\} \right\}}.
$$ 
Let say that our interpolating spline of function $g$ is
$$
s(x):=I[g](x)=\sum_{i=1}^nc_iT_i^k(x),
$$
with $c:=[c_1,\dots,c_n]^{\rm T}=A^{-1}d$, where $d:=[g(\tau_1),\dots,g(\tau_n)]^{\rm T}$. Then
{\setlength\arraycolsep{2pt}
\begin{eqnarray*}
\|I[g]\|_\infty & = & \|\sum_{i=1}^nc_iT_i^k\|_\infty \leqslant \max_{x\in[a,b]}\sum_{i=1}^n|c_i|T_i^k(x) \\
 & \leqslant & \|c\|_\infty\sum_{i=1}^n T_i^k(x) = \|c\|_\infty = \|A^{-1}d\|_\infty \\
 & \leqslant & \|A^{-1}\|_\infty\|d\|_\infty\leqslant \|A^{-1}\|_\infty\|g\|_\infty,
\end{eqnarray*}}%
for all $g\in C([a,b])\backslash \{0\}$, so
$$
\|I\|_\infty \leqslant \|A^{-1}\|_\infty.
$$
Since $\|A\|_\infty=1$, because the B-splines make partition of unity, for the condition number of $A$ we have $\kappa_\infty(A)=\|A^{-1}\|_\infty$, and
$$
\|I\|_\infty \leqslant \kappa_\infty(A).
$$
Like with polynomial splines, $\kappa_\infty(A)$ depends on the choice of the interpolation points $(\tau_i)_{i=1}^n$, but if we choose them like the nodes of the CCC--Schoenberg operator, to fulfil $\tau_i\in[t_{i+1},t_{i+k-1}]$, and not too close to each other, like for an example the CCC--Greville points, we expect $\kappa_\infty(A)$ to be decent.

\section{Collocation and quasi-collocation}

Now we are back to the problem (\ref{boundvalue}), with $f\in C[a,b]$, but we will specify the operator $L_2$ with
\begin{equation} \label{defL2}
L_2=D_2 D_1 D_0
\end{equation}
where $D_0$, $D_1$ and $D_2$ are defined as in (\ref{defDj}) for some $\sigma_2$, $\sigma_3$ and $u_1$, and from now on, we cancel the constraint $u_1\equiv 1$. Let $\cs_2:=\ker{(L_2)}={\rm span}\,{\{u_1,u_2\}}$ with $u_2$ as in (\ref{CCCs}).

For initial value problems, we already have the Green's function, given by Theorem 7 from \cite{bosnerroginaschoenberg}, the same we use to define the B-splines themselves. 
In our case, as for the linear differential operator, we can construct the Green's function also for the boundary value problems (like for tension splines in \cite{marusicroginasharp}).
\begin{theorem} \label{Gaussbv}
The solution to the problem (\ref{boundvalue}), with (\ref{defL2}), is of the form
\begin{equation} \label{solutionGreeen}
y(x)=u(x)+u_1(x)\int_a^b G(x,\tau)f(\tau)\,d\sigma_3(\tau),
\end{equation}
for $x\in[a,b]$, where $u\in\cs_2$ is such that $u(a)=\alpha$ and $u(b)=\beta$, 
and $G$ is the generalized Green's function:
$$
G(x,y)=\left\{ \begin{array}{ll}
\displaystyle -\frac{\int_a^y d\sigma_2 \int_x^b d\sigma_2}{\int_a^b d\sigma_2}, & a\leqslant y\leqslant x\leqslant b \\[1em]
\displaystyle -\frac{\int_a^x d\sigma_2 \int_y^b d\sigma_2}{\int_a^b d\sigma_2}, & a\leqslant x\leqslant y\leqslant b
\end{array} \right. .
$$
\end{theorem}
%
\begin{proof}
The unique function $u$ from (\ref{solutionGreeen}) always exists by Remark 1 from \cite{bosnerroginaschoenberg}. 

First, we want to prove that $y$ satisfies the differential equation. Since 
{\setlength\arraycolsep{2pt}
\begin{eqnarray}
y(x) & = & u(x)-u_1(x)\frac{\int_x^b d\sigma_2}{\int_a^b d\sigma_2} \int_a^x f(\tau)\,d\sigma_3(\tau) \int_a^\tau d\sigma_2 \nonumber \\ 
 & & -u_1(x)\frac{\int_a^x d\sigma_2}{\int_a^b d\sigma_2} \int_x^b f(\tau)\,d\sigma_3(\tau) \int_\tau^b d\sigma_2, \label{yxexpr}
\end{eqnarray}}%
and since the same rule for the generalized derivative of the product of functions is valid, as for the ordinary derivative, we have	
{\setlength\arraycolsep{2pt}
\begin{eqnarray*}
L_2 y(x) & = & \frac{1}{\int_a^b d\sigma_2} D_2\bigg[-\left(D_1\int_x^b d\sigma_2 \right)\cdot \int_a^x f(\tau)d\sigma_3(\tau) \int_a^\tau d\sigma_2 \\
 & & -\int_x^b d\sigma_2 \cdot D_1\left(\int_a^x f(\tau)d\sigma_3(\tau) \int_a^\tau d\sigma_2 \right) \\
 & & -\left(D_1\int_a^x d\sigma_2\right)\cdot\int_x^b f(\tau)\,d\sigma_3(\tau) \int_\tau^b d\sigma_2 \\
 & & -\int_a^x d\sigma_2 \cdot D_1\left(\int_x^b f(\tau)\,d\sigma_3(\tau) \int_\tau^b d\sigma_2 \right) \bigg] \\
 & = & \frac{1}{\int_a^b d\sigma_2} D_2\bigg[\int_a^x f(\tau)d\sigma_3(\tau) \int_a^\tau d\sigma_2 \\
 & & -\int_x^b d\sigma_2 \cdot D_1\left(\int_a^x f(\tau)d\sigma_3(\tau) \int_a^\tau d\sigma_2 \right) -\int_x^b f(\tau)\,d\sigma_3(\tau) \int_\tau^b d\sigma_2 \\
 & &  -\int_a^x d\sigma_2 \cdot D_1\left(\int_x^b f(\tau)\,d\sigma_3(\tau) \int_\tau^b d\sigma_2 \right)\bigg].
\end{eqnarray*}}%
Because by the Mean value theorem for Lebesgue--Stieltjes integrals
$$
D \int_a^x g(\tau)\, d\sigma(\tau) = \lim_{\delta \to 0^+}\frac{\int_x^{x+\delta} g(\tau)\, d\sigma(\tau)}{\int_x^{x+\delta}d\sigma(\tau)} = \lim_{\delta \to 0^+}\frac{g(\xi_x) \int_x^{x+\delta} d\sigma(\tau)}{\int_x^{x+\delta}d\sigma(\tau)} = g(x)
$$
for any measure $\sigma$, any continuous function $g$ for which the above integrals exist, and some $\xi_x\in(x,x+\delta)$, it follows that
{\setlength\arraycolsep{2pt}
\begin{eqnarray}
L_2 y(x) & = & 	\frac{1}{\int_a^b d\sigma_2} \bigg[f(x)\int_a^x d\sigma_2 \\
 & & -\left(D_2\int_x^b d\sigma_2\right)\cdot D_1\left(\int_a^x f(\tau)d\sigma_3(\tau) \int_a^\tau d\sigma_2 \right) \nonumber \\
 & & -\int_x^b d\sigma_2 \cdot D_2D_1\left(\int_a^x f(\tau)d\sigma_3(\tau) \int_a^\tau d\sigma_2 \right)+f(x)\int_x^b d\sigma_2 \nonumber \\
 & &  -\left(D_2\int_a^x d\sigma_2\right) \cdot D_1\left(\int_x^b f(\tau)\,d\sigma_3(\tau) \int_\tau^bd\sigma_2\right) \label{l2yexpr} \\
 & & -\int_a^x d\sigma_2 \cdot D_2D_1\left(\int_x^b f(\tau)\,d\sigma_3(\tau) \int_\tau^bd\sigma_2\right)\bigg]. \nonumber
\end{eqnarray}}%
According to the definition of the generalized derivatives, it is trivial to see that
$$
D_1g(x)=D_2g(x)\cdot D_1\sigma_3(x) \quad {\rm and} \quad D_2g(x)=D_1g(x)\cdot D_2\sigma_2(x),
$$
which yields
{\setlength\arraycolsep{2pt}
\begin{eqnarray}
\left(D_2\int_x^b d\sigma_2\right)\!\!\! & \cdot & \!\!\!D_1\left(\int_a^x f(\tau)d\sigma_3(\tau) \int_a^\tau d\sigma_2 \right) \nonumber \\
 & + & \left(D_2\int_a^x d\sigma_2\right) \cdot D_1\left(\int_x^b f(\tau)\,d\sigma_3(\tau) \int_\tau^b d\sigma_2\right) \nonumber \\
 & = & -D_2\sigma_2(x)\cdot D_2\left(\int_a^x f(\tau)d\sigma_3(\tau) \int_a^\tau d\sigma_2 \right)\cdot D_1\sigma_3(x) \nonumber \\
 & & +D_2\sigma_2(x)\cdot D_2\left(\int_x^b f(\tau)\,d\sigma_3(\tau) \int_\tau^b d\sigma_2\right)\cdot D_1\sigma_3(x) \label{l2y1} \\
 & = & -D_2\sigma_2(x)\cdot f(x)\int_a^x d\sigma_2\cdot D_1\sigma_3(x) \nonumber \\
 & & -D_2\sigma_2(x)\cdot f(x)\int_x^bd\sigma_2\cdot D_1\sigma_3(x) \nonumber \\
 & = & -D_2\sigma_2(x)\cdot D_1\sigma_3(x)\cdot f(x)\int_a^b d\sigma_2, \nonumber
\end{eqnarray}}%
and in the same way it can be shown that
{\setlength\arraycolsep{1pt}
\begin{eqnarray}
\int_x^b d\sigma_2 & \cdot & D_2D_1\left(\int_a^x f(\tau)d\sigma_3(\tau) \int_a^\tau d\sigma_2 \right) + \int_a^x d\sigma_2 \label{l2y2} \\
 & \cdot & D_2D_1\left(\int_x^b f(\tau)\,d\sigma_3(\tau) \int_\tau^bd\sigma_2\right)
 = D_2\sigma_2(x)\cdot D_1\sigma_3(x)\cdot f(x)\int_a^b d\sigma_2. \nonumber
\end{eqnarray}}%
When we put (\ref{l2y1}) and (\ref{l2y2}) back in (\ref{l2yexpr}), we get
$$
L_2y(x)=f(x).
$$
Trivially we can show that $y$ expressed by (\ref{yxexpr}) satisfies the boundary conditions.
\end{proof}
As we mentioned before, the collocation method corresponds to the interpolation problem of the function $f$, so we seek the spline $s\in\cs$ from some spline space $\cs$, such that 
\begin{equation} \label{coloceq}
L_2s=I[f],
\end{equation}
but we offer also an another approach,
to find $s\in\cs$ such that
\begin{equation} \label{quasicoloceq}
L_2s = S[f],
\end{equation}
where both, the interpolation operator and the CCC--Schoenberg operator, are associated with the second reduced system. 
To be precise, we take the measures $d\sigma_2$ and $d\sigma_3$, which define the operator $L_2$, as the first two measures in the measure vector $d\sig$ associated with $\cs$, and we add to them at least one more arbitrary measure. That way, if $\cs(k,d\sig,u_1)$ denotes the underlying CCC--space, the CCC--space associated with the second reduced system $\cs(k-2,d\sig^{(2)})$ is at least two-dimensional. We will discus the possible choice of this added measures a little bit later.

For the spline space, we take the extend partition $\TT$ with all interior knots of multiplicity one, to assure the splines space associated with the second reduced space is at least continuous. Now, our spline space is $\cs:=\cs(k,d\sig,u_1,\TT)$, so the operators $I$ in (\ref{coloceq}) and $S$ in (\ref{quasicoloceq}) are from the spline space $\cs^{(2)}:=\cs(k-2,d\sig^{(2)},\mathbf{1},\TT)$.


%
%

Mostly we use the CCC--Greville points ($\xi_{i,2}^k$) for the CCC--Schoenberg operator, since we know how to calculate them, and although this calculation is complex, for the quasi-collocation we need these Greville points for the spline space of order $k-2$, which is usually not very high. To calculate them, for an example, for order 3 or 4 is often not that complicated. We suggest also to use them as collocation points.

The next question is how good these two approximations are. The answer gives the next theorem.
\begin{theorem} \label{bvquadconv}
Let $f$ from (\ref{boundvalue}), $\es_2$ and $\es_3$ satisfy the conditions of Theorem \ref{thmvdquadr}, only for the spline space $\cs^{(2)}$ instead. Then for the the spline $s\in\cs$ for which (\ref{quasicoloceq}) holds, we have
\begin{equation} \label{bvquadconveq}
\|s-y\|_\infty\leqslant \bar{C}_S \bar{h}_{(d\es_2,d\bar{\sigma}_3,\TT)}^2,
\end{equation}
where $\bar{C}_S$ does not depend on $\TT$. If $f$ is such that $L_{2,i}f$ is continuous on $[a,b]$ for $i=0,1,\dots,k-2$, where the generalized derivatives $L_{2,i}$ are associated with the second reduced measure vector $d\sig^{(2)}=(d\sigma_4,\dots,d\sigma_k)^{\rm T}$, 
then $s\in\cs$ for which (\ref{coloceq}) holds satisfies 
\begin{equation} \label{bvcolocconverg}
\|s-y\|_\infty\leqslant \bar{C}_I \hh_4\cdots\bar{\hh}_{k+1}\|L_{2,k-2}f\|_\infty,
\end{equation}
for an arbitrary measure $d\es_{k+1}$, where $\bar{C}_I$ does not depend on $\TT$.
\end{theorem}
\begin{proof}
Let the operator $Q$ be equal to either $I$ or $S$. 
According to Theorem \ref{Gaussbv}, for both, the exact solution and the approximation, we have
{\setlength\arraycolsep{2pt}
\begin{eqnarray}
y(x) & = & u(x)+u_1(x)\int_a^b G(x,\tau)f(\tau)\,d\sigma_3(\tau) \nonumber \\
s(x) & = & u(x)+u_1(x)\int_a^b G(x,\tau)L_2s(\tau)\,d\sigma_3(\tau) \nonumber \\
 & = & u(x)+u_1(x)\int_a^b G(x,\tau)Q[f](\tau)\,d\sigma_3(\tau) \label{sxGreen}
\end{eqnarray}}%
for some $x\in[a,b]$. Therefore
{\setlength\arraycolsep{2pt}
\begin{eqnarray*}
|s(x)-y(x)| & = & \left|u_1(x)\int_a^b G(x,\tau)\big(Q[f](\tau)-f(\tau)\big)d\sigma_3(\tau)\right| \\
 & \leqslant & |u_1(x)|\,|Q[f](\xi)-f(\xi)|\,\left|\int_a^b G(x,\tau)\,d\sigma_3(\tau)\right| \\
 & \leqslant & \|u_1\|_\infty \|Q[f]-f\|_\infty \left| \int_a^x G(x,\tau)\,d\sigma_3(\tau)+\int_x^b G(x,\tau)\,d\sigma_3(\tau)\right| \\
 & \leqslant & 2\int_a^b d\sigma_2\int_a^b d\sigma_3 \, \|u_1\|_\infty \|Q[f]-f\|_\infty,
\end{eqnarray*}}%
where we used the Mean value theorem for Lebesgue--Stieltjes integrals again. 
Theorem \ref{thmvdquadr}  gives then (\ref{bvquadconveq}), while (\ref{interpolgrocj1}) and Theorem \ref{tmjackson}, the inequality (\ref{bvcolocconverg}).
\end{proof}
\begin{remark} \label{remarkcolqcol}
Actually, we could generalize Theorem \ref{bvquadconv} to any operator $Q$ for which we know the error bounds.
If we use the CCC--Greville points for the Schoenberg operator, the error bound is $O(\hh_4\hh_5)$. 
According to Theorem \ref{tmjackson}, if the function $f$ is less smooth, concerning generalized derivatives, the error bound of the collocation method would be of lower order (``order'' here means the number of factors $\hh_j$), but the same would be true for the quasi-collocation method if the second generalized derivative of $f$ in $\cs_2$ is not at least piecewise continuous.
\end{remark}

Now we can say how to choose the added measures. They should be chosen to reduce the error of function $f$, or its bounds, of the approximation either by interpolation or the Schoenberg operator, as we will see in examples. So the question of the choice of the added measures becomes, more or less, the question how some spline space is good for approximating some given function.
\subsection{Calculation of the approximation} \label{calcaprox}
In this subsection we will offer two ways for calculating the approximation $s$. In the first way, we explicitly calculate the de\thinspace{}Boor points, {\it i.e.} the coefficients of $s$ expressed by the B-splines.

Let 
\begin{equation} \label{sviadbpoints}
s(x)=\sum_{i=1}^n a_i T_i^k(x)
\end{equation}
for some $x\in[a,b]$, then by Theorem \ref{derivativethm}  
{\setlength\arraycolsep{2pt}
\begin{eqnarray}
L_1s(x) & = & \sum_{i=2}^n d_i T_i^{k-1}(x) = \sum_{i=2}^n \frac{a_i-a_{i-1}}{C_i^{k-1}}  T_i^{k-1}(x), \label{l1seq} \\
L_2s(x) & = & \sum_{i=3}^n f_i T_i^{k-2}(x) = \sum_{i=3}^n \frac{d_i-d_{i-1}}{C_i^{k-2}} T_i^{k-2}(x). \label{l2seq}
\end{eqnarray}}%
If we want that $s$ satisfies (\ref{coloceq}) or (\ref{quasicoloceq}), then either we calculate $(f_i)_{i=3}^n$ by solving the interpolation linear system, or $f_i:=f(\zeta_i)$ with $\zeta_i$ from Definition \ref{Schoenbergdef}. Either way, we know the coefficients $f_i$, $i=3,\dots,n$.  Since $f_i=\frac{d_i-d_{i-1}}{C_i^{k-2}}$, it follows that
{\setlength\arraycolsep{2pt}
\begin{eqnarray*}
d_i & = & d_{i-1}+f_i C_i^{k-2} \\
 & = & d_2+\sum_{j=3}^i f_jC_j^{k-2}
\end{eqnarray*}}%
for $i=3,\dots,n$, and in the same way
{\setlength\arraycolsep{2pt}
\begin{eqnarray}
a_i & = & a_{i-1}+d_iC_i^{k-1} \nonumber \\
 & = & a_1+\sum_{\ell=2}^i d_\ell C_\ell^{k-1} \nonumber \\
 & = & a_1+\sum_{\ell=2}^iC_\ell^{k-1}\left(d_2+\sum_{j=3}^\ell f_jC_j^{k-2}\right) \label{aiexp}
\end{eqnarray}}%
for $i=2,\dots,n$. That way, the only unknowns are $a_1$ and $d_2$, which can be found from the boundary values. From $s(a)=\alpha$ we get 
$$
a_1=\frac{\alpha}{u_1(a)},
$$ 
and from $s(b)=\beta$, $a_n=\frac{\beta}{u_1(b)}$. The expression (\ref{aiexp}) gives
$$
a_n = \frac{\alpha}{u_1(a)} +d_2\underbrace{\sum_{\ell=2}^n C_\ell^{k-1}}_{=:A_1}+\underbrace{\sum_{\ell=2}^n C_\ell^{k-1}\sum_{j=3}^n f_jC_j^{k-2}}_{=:A_2} = \frac{\beta}{u_1(b)}
$$
which finally leads to
$$
d_2=\frac{\frac{\beta}{u_1(b)}-\frac{\alpha}{u_1(a)}-A_2}{A_1}.
$$

The alternative way to calculate the coefficients $\{a_i\}_{i=1}^n$ is to formulate the linear system from (\ref{l1seq}) and (\ref{l2seq}):
$$
\frac{\frac{a_i-a_{i-1}}{C_i^{k-1}}-\frac{a_{i-1}-a_{i-2}}{C_{i-1}^{k-1}}}{C_i^{k-2}}=f_i,
$$
for $i=3,\dots n$, together with the boundary conditions. Of course, the collocation coefficients $\{a_i\}_{i=1}^n$ can also be calculated in the usual way, by forming the collocation matrix.

We suggest also another approach, which calculates the value of the spline $s$ (\ref{sviadbpoints}) at some point directly, without explicitly knowing the de\thinspace{}Boor points $a_i$, and which is based on the Green's function and (\ref{sxGreen}). The idea is to find a numerical integration (if it is possible) which can calculate the integrals in (\ref{sxGreen}) (having in mind the expression (\ref{yxexpr})) exactly on each subinterval of the extended partition $\TT$. This way produces a more stable algorithm, than the first one, but also since the calculation of the B-splines, or the evaluation of the spline, can be expensive, that way we only need the B-splines of two order less. The same idea can be applied also to the initial value problems, since as we mentioned before, the solution of it can be represented as in (\ref{solutionGreeen}).

\section{Examples} \label{secexample}

Bellow we illustrate both kind of approximations, and comment several issues mentioned before.

\begin{example}{1} For the first example we chose the singular boundary value problem
{\setlength\arraycolsep{2pt}
\begin{eqnarray}
	D(\sqrt{x}Dy(x)) & = & f(x), \label{examplesingulareq}\\
	y(0)=y(1) & = & 0, \nonumber
\end{eqnarray}}%
with two choices of the function $f$:
\begin{itemize}
\item[a)] 
$f(x)=x^2+x^\frac{5}{2}$,
\item[b)] $\displaystyle f(x)=\frac{e^{100}(1-e^{-100x})}{e^{100}-1}$.
\end{itemize}
For more details and motivation for this problem, see \cite{bosnerroginasingular}. 

So, $u_1\equiv 1$, and the first two measures are: $d\sigma_2(\tau_2)=\frac{d\tau_2}{\sqrt{\tau_2}}$ and $d\sigma_3(\tau_3)=d\tau_3$, and we append them with the same two measures, getting 
$d\sig=(\frac{d\tau_2}{\sqrt{\tau_2}}, d\tau_3,$ $\frac{d\tau_4}{\sqrt{\tau_4}}, d\tau_5)^{\rm T}$. That way, our CCC--space is $\cs(k,d\sig)={\rm span}\,\{1,x^{\frac{1}{2}},x^{\frac{3}{2}},$ $x^2,x^3\}$. The CCC--space associated to the second reduced system is $\cs(k-2,d\sig^{(2)})={\rm span}\,\{1,x^{\frac{1}{2}},x^{\frac{3}{2}}\}$, and all about how to calculate the associated singular splines can be found in \cite{bosnerroginasingular}.


We will calculate the collocation and quasi-collocation approximation via Green's function and Gaussian formula.
In our case, (\ref{sxGreen}) is equal to
{\setlength\arraycolsep{2pt}
\begin{eqnarray}
	s(x) & = & \int_0^1 G(x,t)Q[f](t)\,dt \nonumber \\
	& = & -2 \Big( (1-\sqrt{x})\int_0^x Q[f](t)\sqrt{t}\,dt+\sqrt{x}\int_x^1  Q[f](t)(1-\sqrt{t})\,dt \Big) \label{singintz} \\
	& = & -4 \Big( (1-\sqrt{x})\int_0^{\sqrt{x}} Q[f](z^2)z^2\,dz+\sqrt{x}\int_{\sqrt{x}}^1 Q[f](z^2)(1-z)z\,dz \Big), \nonumber 
\end{eqnarray}}%
by substitution $t=z^2$. Both integrals in the last equation are divided on the integrals over the subintervals, to be precise, if $x\in[t_i,t_{i+1})$, then the first integral is expressed as
$$
\int_0^{\sqrt{x}} = \sum_{j=k}^{i-1} \int_{\sqrt{t_j}}^{\sqrt{t_{j+1}}} + \int_{\sqrt{t_i}}^{\sqrt{x}},
$$
and analogously the second integral. 
Since, after the substitution, the function being integrated by every such integral is polynomial of degree 5, we apply the appropriate Gauss-Legendre integration formula to calculate the integrals in (\ref{singintz}) exactly.

For both, the collocation and quasi-collocation, we use the CCC--Greville points, which are in $\cs^{(2)}$ equal to
$$
\tau_i=\zeta_i=\left(\frac{2}{3}\frac{t_{i+1}+\sqrt{t_{i+1}}\sqrt{t_{i+2}}+t_{i+2}}{\sqrt{t_{i+1}}+\sqrt{t_{i+2}}}\right)^2.
$$

We will consider now each of the proposed functions $f$ separately.

\begin{itemize}
\item[a)] Table \ref{tablgressingular} brings us maximal errors of both approximations for several equidistant partitions, where $n$ is dimension of the underlying spline space. We can clearly see the difference between the approximation by collocation, which is of order 3, in the sense of Theorem \ref{bvquadconv}, and the one by quasi-collocation of order 2. 
We also checked the estimation of the condition number of the interpolation matrix by the Lapack routine \texttt{dgbcon} for all partitions, and they are all around $\kappa_\infty(A)\approx 2.41$. 
\begin{table} 
\begin{center}
\begin{tabular}{|c|c|c|}
\hline
$n$ & quasi-collocation & collocation \\
\hline
     24 &               0.48545E-03 &               0.69229E-06 \\
     44 &               0.12130E-03 &               0.64411E-07 \\
     84 &               0.30320E-04 &               0.59423E-08 \\
    164 &               0.75792E-05 &               0.54331E-09 \\
    324 &               0.18948E-05 &               0.49382E-10 \\
    644 &               0.47368E-06 &               0.44657E-11 \\
   1284 &               0.11842E-06 &               0.40257E-12 \\
   2564 &               0.29605E-07 &               0.36190E-13 \\
   5124 &               0.74012E-08 &               0.32283E-14 \\
  10244 &               0.18503E-08 &               0.54123E-15 \\
  20484 &               0.46258E-09 &               0.45797E-15 \\
  40964 &               0.11565E-09 &               0.10825E-14 \\
  81924 &               0.28911E-10 &               0.10825E-14 \\
 163844 &               0.72278E-11 &               0.11796E-14 \\
\hline
\end{tabular}
\caption{The maximal errors of quasi-collocation and collocation with singular splines for the case a).} \label{tablgressingular}
\end{center}
\end{table}
\item[b)] In this example we will first observe just one equidistant partition, which gives the spline space of dimension $n=15$, and the errors made by both approximations in Figure \ref{slikasinggres}. We can notice that they are, unlikely to the case a), of the same order of magnitude. 
\begin{figure}[h!] 
    \centering
    \includegraphics[width=90mm]{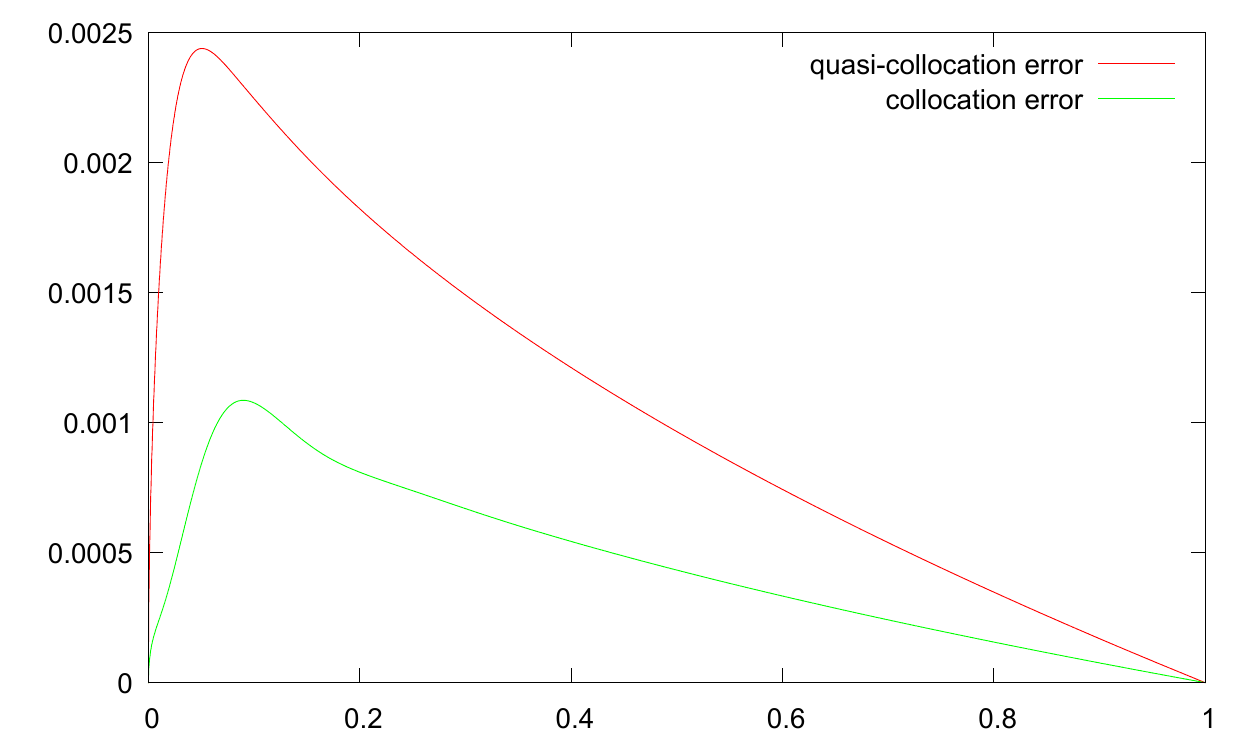}
    \caption{The errors of quasi-collocation and collocation with singular splines for the case b) and $n=15$.} \label{slikasinggres}
\end{figure}
The reason for it is that the second and third generalized derivatives of $f$ are unbounded around zero, so the interpolation of $f$ isn't that good, as we can see in Figure \ref{slikasingf}. The interpolation of $f$ oscillates, while with the approximation by CCC--Schoenberg operator, thanks to its shape preserving properties (see \cite{bosnercrnkovicskificimage,mazurechebschoenberg}), this is not the case.
Therefore, the collocation method is not that much more precise than quasi-collocation, as in case a), what we can see in Table \ref{tablgressingularb} (see Remark \ref{remarkcolqcol}).
\begin{table} 
\begin{center}
\begin{tabular}{|c|c|c|}
\hline
$n$ & quasi-collocation & collocation \\
\hline
     24 &               0.10289E-02 &               0.33555E-03 \\
     44 &               0.39974E-03 &               0.30239E-04 \\
     84 &               0.10698E-03 &               0.18377E-04 \\
    164 &               0.22180E-04 &               0.40525E-05 \\
    324 &               0.41036E-05 &               0.69207E-06 \\
    644 &               0.81481E-06 &               0.11379E-06 \\
   1284 &               0.27664E-06 &               0.19011E-07 \\
   2564 &               0.81434E-07 &               0.32334E-08 \\
   5124 &               0.22471E-07 &               0.56262E-09 \\
  10244 &               0.59872E-08 &               0.97335E-10 \\
  20484 &               0.15618E-08 &               0.17194E-10 \\
  40964 &               0.40189E-09 &               0.30431E-11 \\
  81924 &               0.10250E-09 &               0.53792E-12 \\
 163844 &               0.25981E-10 &               0.95002E-13 \\
\hline
\end{tabular}
\caption{The maximal errors of quasi-collocation and collocation with singular splines for the case b).} \label{tablgressingularb}
\end{center}
\end{table}
For this function $f$ we could maybe make some better choice for the added measures, 
but we wanted to present a case when the collocation method is not much better than quasi-collocation.
\begin{figure}[h!] 
    \centering
    \includegraphics[width=90mm]{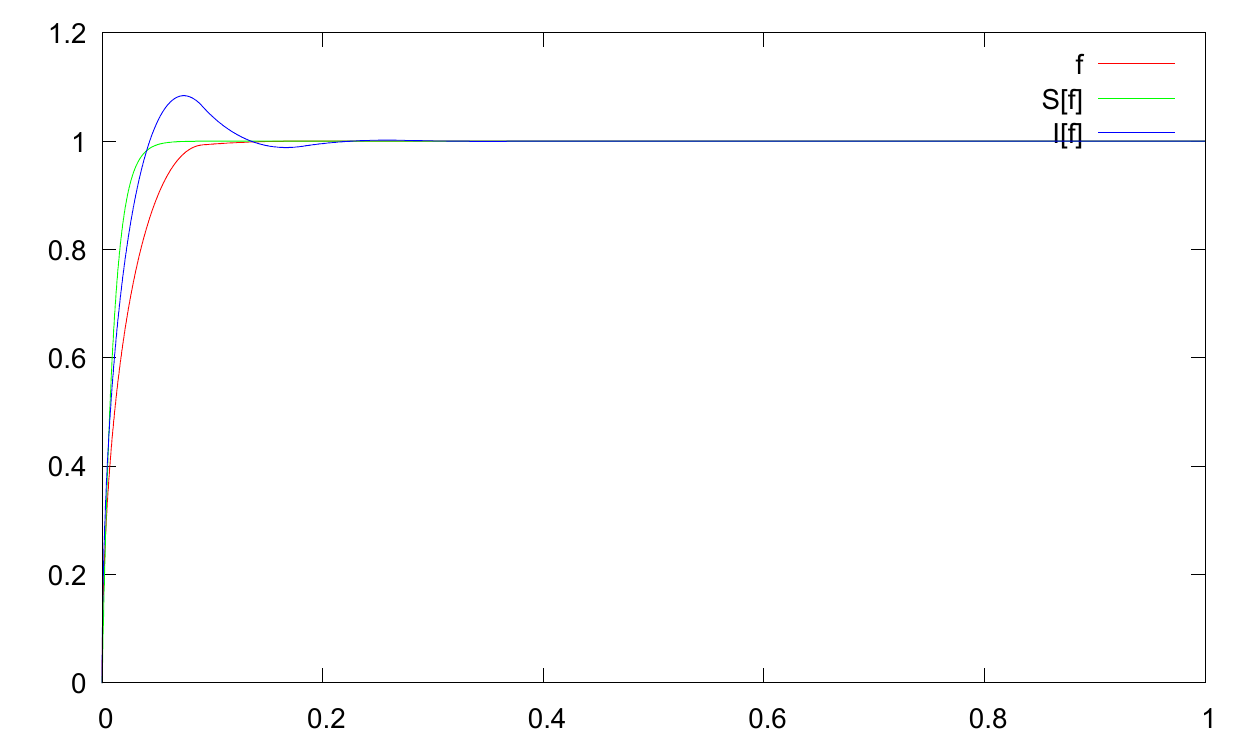}
    \caption{The function $f$ for the case b) and its approximations $S[f]$ and $I[f]$.} \label{slikasingf}
\end{figure}
\end{itemize}

\end{example}

\begin{example}{2}
Next, we solve the boundary value problem
{\setlength\arraycolsep{2pt}
\begin{eqnarray*}
	D^2 y(x)-p^2 y(x) & = & \frac{1}{\cosh{(p x)}} D\left( \cosh{(px)}^2 D\frac{y}{\cosh{(p x)}} \right) = x^2, \\
	y(0)=y(1) & = & 0, \nonumber
\end{eqnarray*}}%
with $p>0$, where we used (\ref{productop}) with $u(x)=\cosh{(p x)}$. This equation for large parameter $p$ is singularly perturbed, giving the solution with one boundary layer, therefore the polynomial splines are, generally, not the proper choice for solving this problem.

Here $u_1(x)=\cosh{(px)}$, and we add one more Lebesgue measure to get the spline space of order 4, with the measure vector $d\sig=(\frac{d\tau_2}{\cosh{(p\tau_2)}^2},\cosh{(p\tau_3)}\,d\tau_3,$ $d\,\tau_4)^{\rm T}$. The algorithms for calculating with such, so called, tension splines, can be found in \cite{bosnerroginatension,bosnerroginanonuniformtension}, although dealing with different CCC--systems, the resulting spline spaces for uniform parameters are the same.

In this example, the space $\cs^{(2)}$ is of order $2$, so the interpolation and the Schoenberg operator give the same approximation ($Q=S=I$). Also, $\cs^{(2)}$ is actually the polynomial spline space, but nevertheless, our main goal is to compare the two ways of calculation of the approximation from Subsection \ref{calcaprox}. So, 
$s_1$ is the approximation of the exact solution $y$, where the de\thinspace{}Boor points are calculated exactly, and $s_2$ is obtained via Green's function and Gaussian formula. 

Now, $s_2$ is
{\setlength\arraycolsep{2pt}
\begin{eqnarray*}
	s_2(x) & = & u_1(x)\int_0^1 G(x,t)Q[f](t)\cosh{(pt)}\,dt \\
	& = & -\frac{1}{p\sinh{p}} \left( \sinh{(p(1-x))}\int_0^x Q[f](t)\sinh{(pt)}\,dt \right. \\ 
	& & \left. +\sinh{(px)}\int_x^1 Q[f](t)\sinh{(p(1-t))}\,dt \right),
\end{eqnarray*}}%
and the integrals over the subintervals can be calculated by using the integral formula that is exact for functions $\cosh{(qx)}$, $x\cosh{(qx)}$, $\sinh{(qx)}$ and $x\sinh{(qx)}$ for some $q>0$, over $[0,1]$, since $Q[f]$ is a linear polynomial spline.
Also, because $\cs^{(2)}$ is the space of polynomial splines, we can numerically verify the order of convergence from Theorem \ref{bvquadconv} in the usual way, denoting them as $p_1$ and $p_2$ for $s_1$ and $s_2$, respectively. 

For $p=10$, the errors and the numerical orders are given in Table \ref{tenserrtblep10}, while for $p=100$ in Table \ref{tenserrtblep100}. The errors for $p=10$ are similar for both calculation algorithms, 
for $p=100$, the first algorithm can not keep pace with the second one, 
and finally for $p=1000$ the first algorithm can not be used at all, so we present only the data for $s_2$ in Table \ref{tenserrtblep1000}.

When the tension splines of order $4$ are concerned, the instability of the first algorithm grows with the growth of the tension parameter $p$, due to the fact that the integrals $C_i^2$ are exponentially growing with $p$, and even more, for the fixed $p$, $C_i^2$ are growing exponentially with $i$, possibly causing overflow, while $C_i^3$ are going smaller with $p$ and $i$ even faster. That way we deal with very big or small values, which vary a lot in their magnitude, obviously causing the instability.
This behaviour is because the density of the measures $d\sigma_2=\frac{d\tau_2}{\cosh{(p\tau_2)}^2}$ and $d\sigma_3=\cosh{(p\tau_3)}\,d\tau_3$, is in one case fast exponentially decaying, while in another case exponentially growing.
This example is chosen to show the advantage of the second algorithm.
\begin{table}
\begin{center}
\begin{tabular}{|c|c|c|c|c|}
\hline
$n$ & $\|s_1-y\|_\infty$ & $p_1$ & $\|s_2-y\|_\infty$ & $p_2$ \\ 
\hline
     23 &     0.412413376460648E-05 &  &     0.412413376454793E-05 & \\
     43 &     0.102848556711096E-05 &   2.003570 &     0.102848556711487E-05 &   2.003570 \\
     83 &     0.256960984222447E-06 &   2.000900 &     0.256960984225049E-06 &   2.000900 \\
     163 &     0.642302149141127E-07 &   2.000225 &     0.642302149119443E-07 &   2.000225 \\
     323 &     0.160569241188825E-07 &   2.000057 &     0.160569241284235E-07 &   2.000057 \\
     643 &     0.401419225810884E-08 &   2.000014 &     0.401419227935920E-08 &   2.000014 \\
     1283 &     0.100354564719005E-08 &   2.000003 &     0.100354565412894E-08 &   2.000003 \\
     2563 &     0.250886219785307E-09 &   2.000001 &     0.250886259683947E-09 &   2.000001 \\
     5123 &     0.627215159150485E-10 &   2.000001 &     0.627215601504971E-10 &   2.000000 \\
     10243 &     0.156804530297705E-10 &   1.999993 &     0.156804100953645E-10 &   1.999998 \\
     20483 &     0.392006956409507E-11 &   2.000016 &     0.392012854469326E-11 &   1.999990 \\
     40963 &     0.980084503138245E-12 &   1.999901 &     0.980073661116521E-12 &   1.999939 \\
     81923 &     0.244914765551441E-12 &   2.000626 &     0.245122932368558E-12 &   1.999385 \\
     163843 &     0.612101515307106E-13 &   2.000437 &     0.613597714305136E-13 &   1.998140 \\
\hline
\end{tabular}
\end{center}
\caption{The errors and numerical orders of $s_1$ and $s_2$ for tension splines of order $4$ and $p=10$.} \label{tenserrtblep10}
\end{table}
\begin{table}
\begin{center}
\begin{tabular}{|c|c|c|c|c|}
\hline
$n$ & $\|s_1-y\|_\infty$ & $p_1$ & $\|s_2-y\|_\infty$ & $p_2$ \\ 
\hline
     23 &     0.903501233085530E-04 &  &     0.507641834927993E-07 & \\
     43 &     0.903501233085530E-04 &   0.000000 &     0.112310617536170E-07 &   2.176317 \\
     83 &     0.444711676796353E+20 & -78.703618 &     0.266114078252636E-08 &   2.077378 \\
     163 &     0.903501233085530E-04 &  78.703618 &     0.654711583612886E-09 &   2.023113 \\
     323 &     0.903501233085530E-04 &   0.000000 &     0.162991554037177E-09 &   2.006062 \\
     643 &     0.903501233085530E-04 &   0.000000 &     0.407045782035943E-10 &   2.001534 \\
     1283 &     0.903501233085530E-04 &   0.000000 &     0.101734313112450E-10 &   2.000385 \\
     2563 &     0.903501233085530E-04 &   0.000000 &     0.254318830263719E-11 &   2.000096 \\
     5123 &     0.903501233085530E-04 &   0.000000 &     0.635786697811627E-12 &   2.000024 \\
     10243 &     0.903501233085530E-04 &   0.000000 &     0.158946640571589E-12 &   2.000000 \\
     20483 &     0.903501233085530E-04 &   0.000000 &     0.397372971116736E-13 &   1.999977 \\
     40963 &     0.444711676796353E+20 & -78.703618 &     0.993564226118432E-14 &   1.999809 \\
     81923 &     0.903501233085530E-04 &  78.703618 &     0.248550637536871E-14 &   1.999073 \\
     163843 &     0.903501233085530E-04 &   0.000000 &     0.623240066326136E-15 &   1.995680 \\
\hline
\end{tabular}
\end{center}
\caption{The errors and numerical orders of $s_1$ and $s_2$ for tension splines of order $4$ and $p=100$.} \label{tenserrtblep100}
\end{table}
\begin{table}
\begin{center}
\begin{tabular}{|c|c|c|}
\hline
$n$ & $\|s_2-y\|_\infty$ & $p_2$ \\ 
\hline
     23 &     0.623000000005656E-09 & \\
     43 &     0.154240187269918E-09 &   2.014053 \\
     83 &     0.371083218788035E-10 &   2.055364 \\
     163 &     0.831544592344676E-11 &   2.157877 \\
     323 &     0.181163400130717E-11 &   2.198502 \\
     643 &     0.420507116972695E-12 &   2.107089 \\
     1283 &     0.102616142194353E-12 &   2.034873 \\
     2563 &     0.254876882039150E-13 &   2.009385 \\
     5123 &     0.636138175827156E-14 &   2.002388 \\
     10243 &     0.158969110661614E-14 &   2.000594 \\
     20483 &     0.397412082863075E-15 &   2.000039 \\
     40963 &     0.993848109210704E-16 &   1.999538 \\
     81923 &     0.248615816722162E-16 &   1.999107 \\
     163843 &     0.622611567879273E-17 &   1.997514 \\
\hline
\end{tabular}
\end{center}
\caption{The errors and numerical orders of $s_2$ for tension splines of order $4$ and $p=1000$.} \label{tenserrtblep1000}
\end{table}
\end{example}

\begin{example}{3}
Our last example is a very simple one:
{\setlength\arraycolsep{2pt}
\begin{eqnarray*}
	D^2 y(x) & = & \frac{\sinh{(px)}}{\sinh{p}}px, \\
	y(0) & = & y(1) = 0, \nonumber
\end{eqnarray*}}%
with $p>0$, which is also a singularly perturbed equation, 
and for a large parameter $p$ the solution, as well as the right hand size $f$, has one boundary layer. In this case, $u_1\equiv 1$, and both measures, $d\sigma_2$ and $d\sigma_3$ are Lebesgue measures. The next question is how to chose the remaining measures. Because of the boundary layer of the solution, the simplest choice of just adding some more Lebesgue measures is not the best one, since, as we mentioned before, 
polynomial splines do not catch that layer well (especially if the number of knots isn't too big). So, we take the following measure vector
$$
d\sig=\left( d\tau_2,d\tau_3,\cosh{(p\tau_4)}\,d\tau_4,\frac{d\tau_5}{\cosh^2(p\tau_5)}\right)^{\rm T},
$$
which gives us the tension splines of order $5$ (piecewisely spanned by $\{1,x,x^2, \cosh{(px)},$ $\sinh{(px)}\}$), which is a reasonable choice, concerning the right hand side of the equation, and the boundary layer. We use 
the second approximation by integrals $s_2$ for quasi-collocation, because for that algorithm we only need to calculate the B-splines of order $3$ from the space $\cs^{(2)}$. 
In this case the approximations by quasi-collocation and collocation are almost the same, because $\cs^{(2)}$ is almost the space of polynomial splines of order 1, so we will skip the collocation.
Then, we compare $s_2$ with the polynomial splines of order 5 $s_{sp}$ and $s_{ip}$, obtained by the first algorithm. The approximation $s_{sp}$ is achieved by quasi-collocation, and $s_{ip}$ by collocation.


Here
{\setlength\arraycolsep{2pt}
\begin{eqnarray}
	s_2(x) & = & \int_0^1 G(x,t)S[f](t)\,dt \nonumber \\
	& = & -  (1-x)\int_0^x S[f](t)\,dt -x \int_x^1 S[f](t)(1-t)\,dt, \label{s2tens5int}
\end{eqnarray}}%
and $S[f]$ is piecewisely spanned by $\{1,\cosh{(px)},\sinh{(px)}\}$. It is possible to find stable integration formulas to calculate integrals in (\ref{s2tens5int}) exactly. For the points $(\zeta_i)_{i=1}^n$ which define the CCC--Schoenberg operator in $\cs^{(2)}$, we chose the Greville abscis\ae{}, where $u_{2,2}(x)=\frac{\sinh{px}}{p}$, so
$$
\zeta_i=\frac{1}{p}\arcsinh{\left(\frac{\sinh{\frac{p(t_{i+1}+t_{i+2})}{2}}}{\cosh{\frac{p(t_{i+2}-t_{i+1})}{2}}}\right)}.
$$

We compare these three approximation for $p=10000$ in Table \ref{tens5errtblep10000}. 
The stagnation in errors of $s_2$ for smaller $n$ is most probably due a specific shape of the solution, which looks almost as a linear function, with a very steep boundary layer. The maximum of the exact solution is around $10^{-5}$, so we see from that table, that the polynomial splines only for large $n$ come a little bit closer to the solution, and even interpolation of $f$, which is of order 3, doesn't help a lot. Also, in the same table we see the connection between large error of $f$ and large error of $y$.
\begin{table} 
\begin{center}
\begin{tabular}{|c|c|c|c|c|} 
\hline
$n$ & $\|s_2-y\|_\infty$ & $\|L_2s_2-f\|_\infty$ & $\|s_{sp}-y\|_\infty$ & $\|L_2s_{sp}-f\|_\infty$ \\ 
\hline 
     24 &       0.200E-07 &       0.368E+00 &       0.200E+01 &       0.974E+04 \\     
     44 &       0.200E-07 &       0.368E+00 &       0.510E+00 &       0.954E+04 \\       
     84 &       0.200E-07 &       0.368E+00 &       0.129E+00 &       0.919E+04 \\       
    164 &       0.200E-07 &       0.368E+00 &       0.323E-01 &       0.861E+04 \\       
    324 &       0.200E-07 &       0.368E+00 &       0.801E-02 &       0.768E+04 \\       
    644 &       0.200E-07 &       0.368E+00 &       0.194E-02 &       0.628E+04 \\       
   1284 &       0.182E-07 &       0.364E+00 &       0.469E-03 &       0.437E+04 \\       
   2564 &       0.890E-08 &       0.287E+00 &       0.142E-03 &       0.246E+04 \\       
   5124 &       0.324E-08 &       0.133E+00 &       0.431E-04 &       0.117E+04 \\       
  10244 &       0.105E-08 &       0.478E-01 &       0.116E-04 &       0.501E+03 \\       
  20484 &       0.288E-09 &       0.175E-01 &       0.296E-05 &       0.186E+03 \\      
  40964 &       0.738E-10 &       0.564E-02 &       0.743E-06 &       0.578E+02 \\       
  81924 &       0.186E-10 &       0.151E-02 &       0.186E-06 &       0.152E+02 \\       
 163844 &       0.465E-11 &       0.384E-03 &       0.465E-07 &       0.383E+01 \\       
\hline
\end{tabular} \\[1em]

\begin{tabular}{|c|c|c|}
\hline
$n$ & $\|s_{ip}-y\|_\infty$ & $\|L_2s_{ip}-f\|_\infty$ \\
\hline 
24 & 0.120E+01 &       0.965E+04 \\
44 & 0.306E+00 &       0.938E+04 \\
84 & 0.776E-01 &       0.893E+04 \\
164 & 0.196E-01 &       0.819E+04 \\
324 & 0.499E-02 &       0.702E+04 \\
644 & 0.132E-02 &       0.531E+04 \\
1284 & 0.340E-03 &       0.312E+04 \\
2564 & 0.547E-04 &       0.121E+04 \\
5124 & 0.452E-05 &       0.333E+03 \\
10244 & 0.318E-06 &       0.610E+02 \\
20484 & 0.205E-07 &       0.873E+01 \\
40964 & 0.129E-08 &       0.113E+01 \\
81924 & 0.808E-10 &       0.864E-01 \\
163844 & 0.505E-11 &       0.135E-01 \\
\hline
\end{tabular}
\end{center}
\caption{The errors of $s_s$, $s_{sp}$ and $s_{ip}$, together with the errors of the approximations of $f$ for $p=10000$.} \label{tens5errtblep10000}
\end{table}    

To illustrate such big differences, we concentrate more closely on the case when $n=24$. First, we observe the approximations of the function $f$, {\it i.e.} $S[f]$ from the spaces of tension and polynomial splines, and $I[f]$ from polynomial spline space. As we can see on Figure \ref{slikatensf}, the tension spline approximation is almost exactly above the function $f$, which is almost everywhere equal to zero, except that very steep boundary layer on the right, while the polynomial ones are quite far away, and the interpolation is oscillating a lot again. 
\begin{figure}[h!] \label{slikatensf}
    \centering
    \includegraphics[width=90mm]{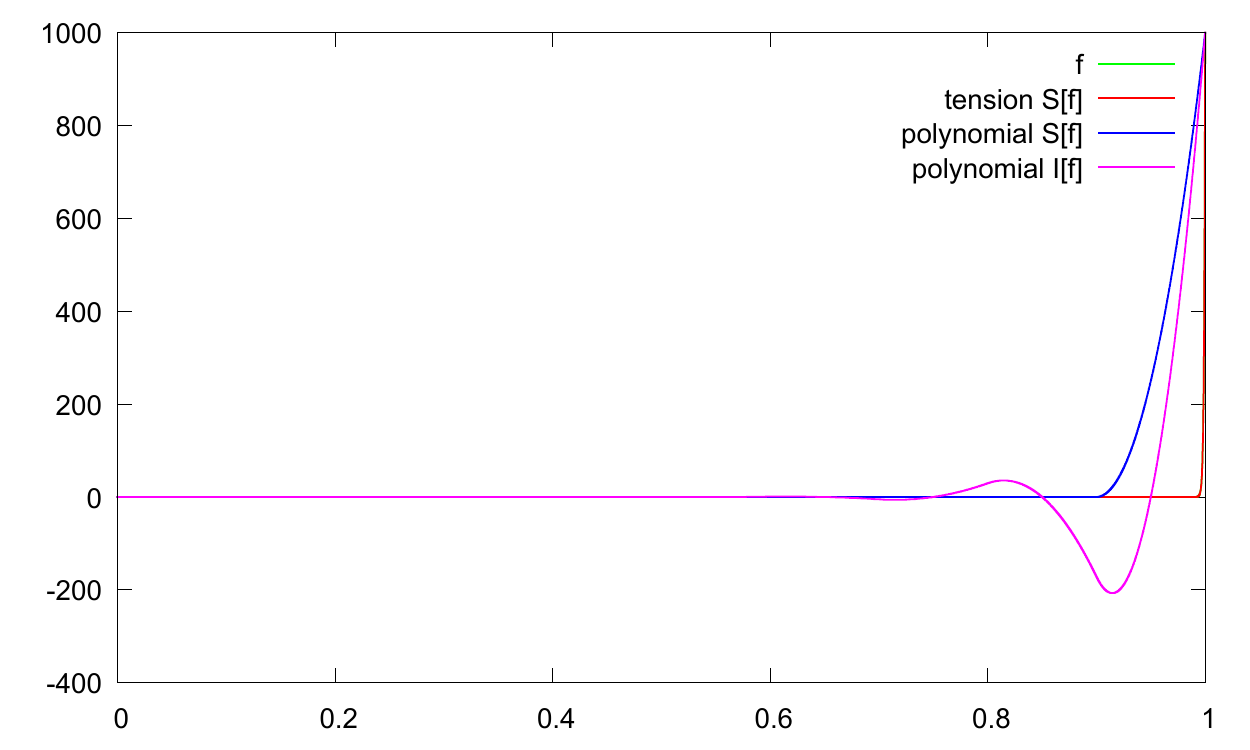}
    \caption{The function $f$ and its approximations.}
\end{figure}
For the approximation of the solution, the situation is even worse, the polynomial splines are completely wrong and way out of scale (the approximation by polynomial quasi-collocation is at least closer in shape), as we can see on Figures \ref{slikatensu1} and \ref{slikatensu2}.
\begin{figure}[h!] 
    \centering
    \includegraphics[width=90mm]{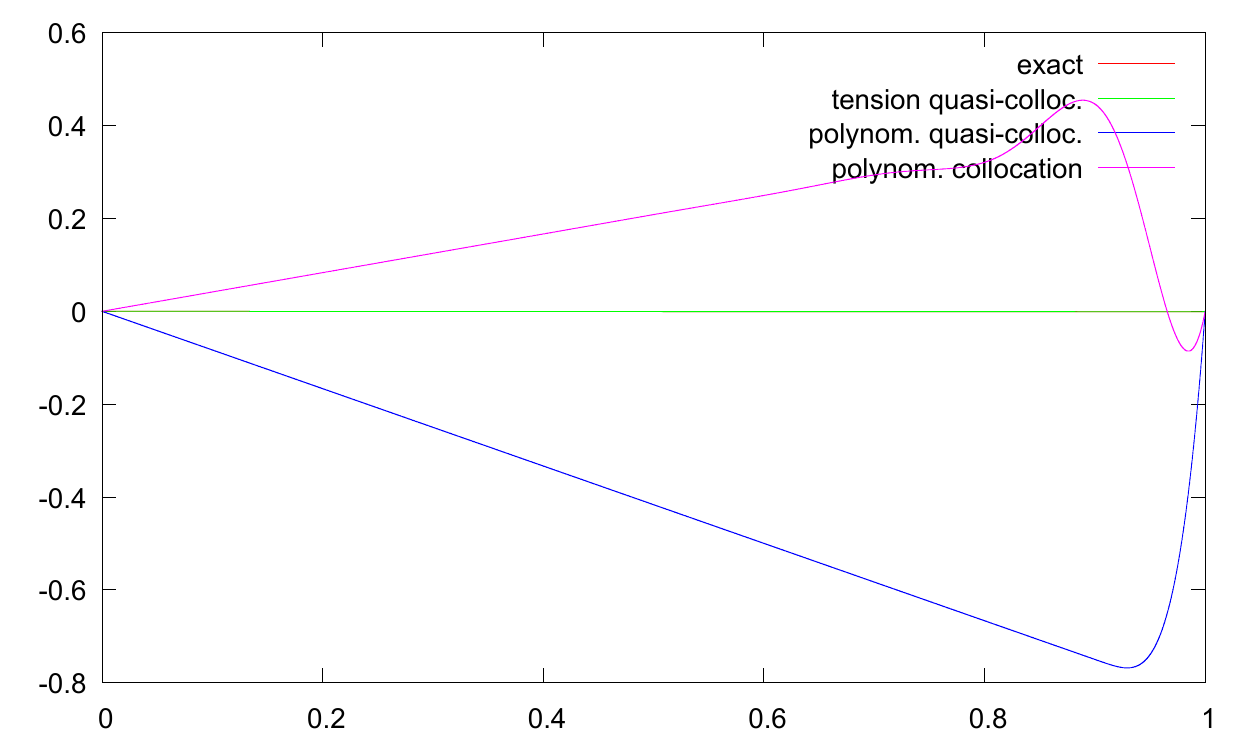}
    \caption{The solution and all three approximations.} \label{slikatensu1}
\end{figure}
\begin{figure}[h!] 
    \centering
    \includegraphics[width=90mm]{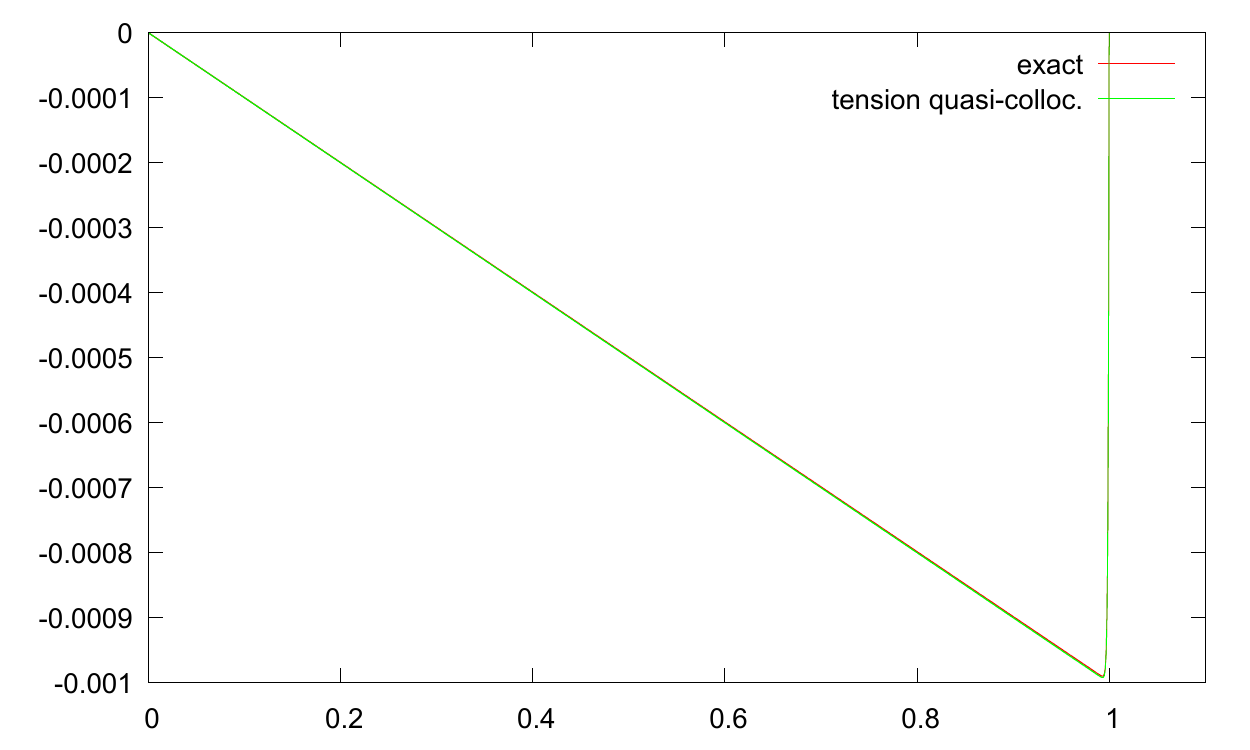}
    \caption{The solution and the tension spline approximation.} \label{slikatensu2}
\end{figure}
%


Therefore, we can conclude, that the simplest choice for added measures, {\it i.e.} the Lebesgue measures, doesn't have to be always good, and as we mentioned before, 
the better approximation of $f$ suggests that the chosen measures will lead to the smaller error of the final solution.

\end{example}

\section{Conclusion}

In this paper we demonstrated one more possible use of the CCC--Schoenberg operators, in solving boundary value problems, with emphasize on the singular and singularly perturbed ODE, and compare it with the collocation method. We proposed the algorithm for calculating both of them, based on the integral of the Green's function, which uses only the B-splines of two orders less. They are, of course, simpler to calculate than the B-splines from the original spline space.
For the mostly used spline spaces, it should not be a problem to find integration formulas, which allow us to calculate the required integrals exactly. When we can numerically stably calculate the B-splines from the original space, the first algorithm based on calculation of de\thinspace{}Boor points, can also be used without a problem, in cases where the integrals $C_i^{k-1}$ and $C_i^{k-2}$ in (\ref{l1seq}) and (\ref{l2seq}) are ``decent''. 
We also show in our examples the cases when the quasi-collocation can almost be as good  as the collocation method. 
Therefore, in our opinion, this quasi-collocation method could be a useful tool for solving these specific problems.


%
%

\bibliographystyle{spmpsci}      
\bibliography{references}   


\end{document}